\documentclass[a4paper,12pt,article]{aplimattpl}

\usepackage[T1]{fontenc}
\usepackage[utf8]{inputenc}

\usepackage{amsmath}
\usepackage{amssymb}
\usepackage{times}

\DeclareMathOperator{\e}{e}

\begin{document}

\pagestyle{aplimat}

\authorsandtitle{REBENDA Josef/CZ}{An application of Bell polynomials in numerical solving of nonlinear differential equations}{1.}

\abstract{Partial ordinary Bell polynomials are used to formulate and prove a version of the Fa\`{a} di Bruno's formula which is convenient for handling nonlinear terms in the differential transformation. Applicability of the result is shown in two examples of solving the initial value problem for differential equations which are nonlinear with respect to the dependent variable.}



\newsavebox{\authors}
\savebox{\authors}{%
\parbox{1.0\textwidth}{%
		\setlength{\parskip}{6pt}%
       
    \textbf{Rebenda Josef, Mgr., PhD.}\\
    Department of Mathematics\\
    Faculty of Electrical Engineering and Communication\\
    Brno University of Technology\\
    Technick\'{a} 8, 616 00 Brno, Czech Republic\\
    E-mail: josef.rebenda@ceitec.vutbr.cz
    }
}

 \begin{aplart}

\section{Introduction}

The ability to find a numerical approximation of solution of a differential equation is important in particular practical applications where it is difficult or even impossible to find analytical solution of the given problem. Plenty of well-established and verified numerical methods for various kinds of problems involving differential equations can be found for instance in monographs \cite{bellen}, \cite{butcher}, \cite{hairer1}, \cite{hairer2}.

Among others, semi-analytical methods convenient for solving differential equations are in the forefront of study in the last two decades. However, the calculations and results are often expressed in a complicated way. We propose an easily applicable approach in this paper.

The differential transformation is closely related to Taylor expansion of real analytic functions with applications to different types of problems of solving differential equations. To indicate recent development in the field we mention several papers from the last three years, e.g. \cite{rebenda1}, \cite{rebenda2}, \cite{rebenda3}, \cite{samajova}, \cite{yang}.

The paper is organized as follows. First we recall basic definitions and formulas of the differential transformation and give a brief overview of Bell polynomials and the Fa\`{a} di Bruno's formula in Section \ref{difftransform}. In Section \ref{main} we develop the theory and prove the main result. Application of the results is shown in Section \ref{examples}.

Convergence, error estimates and stability of Taylor series based methods is thoroughly discussed in literature on numerical methods, therefore we do not include such topics in the paper. However, interesting results on a-priori error bounds are published in recent paper \cite{warne} and the references cited therein.

\section{Preliminaries}
\label{difftransform}

In this section we recall basic definitions and formulas of the differential transformation as well as notions and results related to partial Bell polynomials.

\textbf{Definition 1}
The differential transformation of a real function $u(t)$ at a point $t_0 \in \mathbb R$ is $\mathcal D \{ u(t) \} [t_0] = \{ U(k) [t_0] \}_{k=0}^{\infty}$,
where $U(k) [t_0]$, the differential transformation of the $k-$th derivative of the function $u(t)$ at $t_0$, is defined as
\begin{equation}\label{2}
U(k) [t_0] = \frac{1}{k!} \left[ \frac{d^ku(t)}{dt^k} \right]_{t=t_0},
\end{equation}
provided that the original function $u(t)$ is analytic in some neighbourhood of $t_0$.

\textbf{Definition 2}
 The inverse differential transformation of $\{ U(k) [t_0] \}_{k=0}^{\infty}$ at $t_0$ is defined as
\begin{equation}\label{3}
u(t) = \mathcal D^{-1} \Bigl\{ \ \{ U(k) [t_0] \}_{k=0}^{\infty} \ \Bigr\} [t_0]= \sum_{k=0}^{\infty}U(k) [t_0] (t-t_0)^k. 
\end{equation}

In real applications the function $u(t)$ is expressed by a finite sum
\begin{equation}\label{3f}
u(t) = \sum_{k=0}^{N}U(k) [t_0] (t-t_0)^k. 
\end{equation}
Plenty of transformation formulas can be derived from Definitions 1 and 2. We recall the following relations which will be used later in illustratory examples.

\textbf{Lemma 1}\label{lemma1}
Assume that $\{ F(k) \}_{k=0}^{\infty}$, $\{ G(k) \}_{k=0}^{\infty}$, $\{ H(k) \}_{k=0}^{\infty}$ and $\{ U_i(k) \}_{k=0}^{\infty}$, $i=1,\dots,m$, are differential transformations of analytic functions $f(t)$, $g(t)$, $h(t)$ and
$u_i(t)$, $i=1,\dots,m$, at $t_0 \in \mathbb{R}$, respectively. Let $n \in \mathbb{N}$ and $\lambda \in \mathbb{R}$. Then
\begin{gather*}
\begin{array}{lllcl}
i) & \text{If} & f(t) = {\displaystyle \frac{d^ng(t)}{dt^n}}, & \text{then} & F(k) = {\displaystyle \frac{(k+n)!}{k!}}G(k+n). \\[3mm]
ii)& \text{If} & f(t) = g(t) \ h(t), & \text{then} & F(k) = \sum\limits_{l=0}^k G(l) \ H(k-l). \\[4mm]
iii)& \text{If} & f(t) = t^n, & \text{then} & F(k) =\delta (k-n), t_0=0, \ \text{where} \ \delta \ \text{is\ the\ Kronecker}\\[2mm]
 \ & \ & \ & \ & \text{delta symbol}. \\[2mm]
iv)& \text{If} & f(t) = e^{\lambda t}, & \text{then} & F(k) = {\displaystyle \frac{\e^{a \lambda} \lambda^k}{k!}}, t_0=a. \\[3mm]
v)& \text{If} & f(t) = \ln t, & \text{then} & F(0) = 0 \ \text{and } \ F(k) = {\displaystyle \frac{(-1)^{k+1}}{k}}  \ \text{for} \ k \geq 1, t_0=1. \\[3mm]
vi)& \text{If} & f(t) = (1+t)^{\lambda}, & \text{then} & F(k) = {\displaystyle \binom{\lambda}{k} = \frac{\lambda (\lambda-1) \ldots (\lambda-k+1)}{k!}}, t_0=0. \\[2mm]
vii) &\text{If} & f(t) = \prod\limits_{i=1}^m u_i(t),& \text{then}&
\end{array}
\\
F(k)= \sum_{s_1=0}^k \ \sum_{s_2=0}^{k-s_1} \ \dots \ \sum_{s_{m-1}=0}^{k-s_1-\ldots - s_{m-2}} U_1(s_1) \dots U_{m-1}(s_{m-1}) \ U_m(k-s_1-\ldots -s_{m-1}).
\end{gather*}

The main disadvantage of most papers with applications of the differential transformation is that there is lack of direct applications on equations with nonlinear terms containing unknown function $u(t)$, e.g. $f(u) = \sqrt{1+u^2}$ or $f(u)= e^{\sin{u}}$. Usually, the "nonlinearity" is represented by terms $u^n$, $n \in \mathbb{N}$, and formula $vii)$ in Lemma 1 is used to transform such terms.

In paper \cite{smardaz} the differential transformation of components containing nonlinear terms is calculated using the so-called Adomian polynomials $A_n$ in which each solution component $u_i$ is replaced  by the corresponding  differential transformation component $U(i), \ i=0,1,2, \dots$. The formula for the differential transformation $F(k)$ of a nonlinear term $f(u)$ is

\begin{align*}\label{PP4}
F(k) &=  \sum_{n=0}^{\infty} A_n(U(0),U(1),\dots, U(n))\delta(k-n) 
 = A_k(U(0),U(1),\dots, U(k))\\
&=  \frac{1}{k!} \frac{d^k}{dt^k}  \left[ f\left(  \sum_{i=0}^{\infty} U(i) t^i \right)\right]_{t=0}, \quad k \geq 0. 
\end{align*}

The first four terms are:
\begin{eqnarray*}\label{A2}
F(0) &=& f(U(0)),  \\
F(1) &=& U(1)f'(U(0)), \\
F(2) &=& U(2)f'(U(0)) +\frac{1}{2!} U^2(1)f''(U(0))\\
F(3) &=& U(3)f'(U(0))+U(1)U(2)f''(U(0))+ \frac{1}{3!} U^3(1) f'''(U(0)).
\end{eqnarray*}

As we can observe, the formula with the so-called Adomian polynomials is in fact the well-known Fa\`{a} di Bruno's formula generalizing the chain rule to higher derivatives. However, there are derivatives of the function $f$ contained in the formula, which means that symbolic derivatives of $f$ need to be calculated and then evaluated. In such situation, one of the big advantages of the differential transformation is lost.

Fortunately, the differential transformation of components containing nonlinear terms can be easily found without calculating and evaluating symbolic derivatives. We will utilize a slightly modified Fa\`{a} di Bruno's formula with not exponential but ordinary Bell polynomials.

For this purpose we recall some necessary notions and results in combinatorics. The proofs are omitted since they can be found in the cited literature \cite{chara} and \cite{comtet}.

\textbf{Definition 3}[\cite{comtet}, p. 133]
\label{pebp}
The partial exponential Bell polynomials are the polynomials\\
$B_{k,l} (x_1, \ldots, x_{k-l+1})$ in an infinite number of variables $x_1, x_2, \ldots$, defined by the series expansion
\begin{equation}\label{pebp1}
\sum\limits_{k \geq l} B_{k,l} (x_1, \ldots, x_{k-l+1}) \frac{t^k}{k!} = \frac{1}{l!} \left( \sum\limits_{m \geq 1} x_m \frac{t^m}{m!} \right)^l, \quad l=0,1,2, \ldots
\end{equation}

\textbf{Lemma 2}[\cite{comtet}, p. 134]
\label{pebp2}
The partial exponential Bell poynomials have integer coefficients, are homogeneous of degree $l$ and weight $k$, and their exact expression is:
\begin{equation}\label{pebp3}
B_{k,l} (x_1, \ldots, x_{k-l+1}) = \sum \frac{k!}{j_1 ! j_2 ! \cdots j_{k-l+1} !} \left( \frac{x_1}{1!} \right)^{j_1} \left( \frac{x_2}{2!} \right)^{j_2} \cdots \left( \frac{x_{k-l+1}}{(k-l+1)!} \right)^{j_{k-l+1}},
\end{equation}
where the summation takes place over all sequences $j_1$, $j_2$, \dots, $j_{k-l+1}$ of non-negative integers such that
\begin{align}
j_1 &+ j_2 + \ldots + j_{k-l+1} = l,\\
j_1 &+ 2 j_2 + \ldots + (k-l+1) j_{k-l+1} = k.
\end{align}

\textbf{Lemma 3}[\cite{chara}, p. 415]
The partial exponential Bell polynomials $B_{k,l} (x_1, \ldots, x_{k-l+1})$,\\
$l=1,2, \ldots$, $k \geq l$,  satisfy the recurrence relation
\begin{equation}\label{pebp4}
B_{k,l} (x_1, \ldots, x_{k-l+1}) = \sum\limits_{i=1}^{k-l+1} \binom{k-1}{i-1} x_i B_{k-i,l-1} (x_1, \ldots, x_{k-i-l+2}),
\end{equation}
where $B_{0,0} = 1$ and $B_{k,0}=0$ for $k \geq 1$.

\textbf{Theorem 1}[\cite{comtet}, p. 138-139]
\label{theorem1}
Let two functions $f(y)$ and $g(x)$ of a real variable be given, $g(x)$ of class $C^{\infty}$ at $x=t_0$, and $f(y)$ of class $C^{\infty}$ at $y=s_0 = g(t_0)$, and let $h(x) = (f \circ g) (x) = f(g(x))$. If we put $g_0 = g(t_0)$, $f_0 = f(s_0) = h_0 = h(t_0) = f(g(t_0))$, $\displaystyle{g_m = \bigg. \frac{d^m g}{dx^m} \bigg\vert_{x=t_0}}$, $\displaystyle{f_l = \bigg. \frac{d^l f}{dy^l} \bigg\vert_{y=s_0}}$, $\displaystyle{h_k = \bigg. \frac{d^k h}{dx^k} \bigg\vert_{x=t_0}}$, then the $k$-th order derivative of $h$ at $t=t_0$ for $k \geq 1$ equals
\begin{equation}\label{fdb}
h_k = \bigg. \frac{d^k h}{dx^k} \bigg\vert_{x=t_0} = \sum\limits_{l=1}^{k} f_l B_{k,l} (g_1, g_2, \ldots, g_{k-l+1}),
\end{equation}
where $B_{k,l}$ are explicitly given by \eqref{pebp3}.

\textbf{Definition 4}[\cite{comtet}, p. 136]
\label{pobp}
The partial ordinary Bell polynomials are the polynomials\\
$\hat{B}_{k,l} (\hat{x}_1, \ldots, \hat{x}_{k-l+1})$ in an infinite number of variables $\hat{x}_1, \hat{x}_2, \ldots$, defined by the series expansion
\begin{equation}\label{pobp1}
\sum\limits_{k \geq l} \hat{B}_{k,l} (\hat{x}_1, \ldots, \hat{x}_{k-l+1}) t^k =\left( \sum\limits_{m \geq 1} \hat{x}_m t^m \right)^l, \quad l=0,1,2, \ldots
\end{equation}

\section{Results}\label{main}

Before we formulate the main theorem, we prove several auxiliary results. The following lemma is crucial in the proof of the main result.

\textbf{Lemma 4}
\label{lemma4}
The relation between the partial exponential Bell polynomials $B_{k,l}$ and the partial ordinary Bell polynomials $\hat{B}_{k,l}$ is
\begin{equation}\label{pobp2}
B_{k,l} (x_1, \ldots, x_{k-l+1}) = \frac{k!}{l!} \hat{B}_{k,l} \left( \frac{x_1}{1!}, \frac{x_2}{2!}, \ldots, \frac{x_{k-l+1}}{(k-l+1)!} \right).
\end{equation}

\textbf{Proof.}
If we denote ${\displaystyle \hat{x}_i = \frac{x_i}{i!}}$ for all $i = 0, 1, 2, \ldots$, and substitute in \eqref{pobp1}, we obtain
\begin{equation}
\sum\limits_{k \geq l} \hat{B}_{k,l} \left( \frac{x_1}{1!}, \frac{x_2}{2!}, \ldots, \frac{x_{k-l+1}}{(k-l+1)!} \right) t^k = \left( \sum\limits_{m \geq 1} \frac{x_m t^m}{m!} \right)^l, \quad l=0,1,2, \ldots
\end{equation}
Multiplying both sides by ${\displaystyle \frac{1}{l!}}$ and each summand on the left side by ${\displaystyle \frac{k!}{k!}}$, we get
\begin{equation}\label{pobp3}
\sum\limits_{k \geq l} \frac{k!}{l!} \hat{B}_{k,l} \left( \frac{x_1}{1!}, \frac{x_2}{2!}, \ldots, \frac{x_{k-l+1}}{(k-l+1)!} \right) \frac{t^k}{k!} = \frac{1}{l!} \left( \sum\limits_{m \geq 1} x_m \frac{t^m}{m!} \right)^l, \quad l=0,1,2, \ldots
\end{equation}
Equating the coefficients of $t^k$ in \eqref{pobp3} and \eqref{pebp1} we get the formula \eqref{pobp2}.

In calculation of the partial ordinary Bell polynomials, the following lemma can be useful.

\textbf{Lemma 5}
The partial ordinary Bell polynomials $\hat{B}_{k,l} (\hat{x}_1, \ldots, \hat{x}_{k-l+1})$,
$l=1,2, \ldots$, $k \geq l$, satisfy the recurrence relation
\begin{equation}\label{pobp4}
\hat{B}_{k,l} (\hat{x}_1, \ldots, \hat{x}_{k-l+1}) = \sum\limits_{i=1}^{k-l+1} \frac{i \cdot l}{k} \hat{x}_i \hat{B}_{k-i,l-1} (\hat{x}_1, \ldots, \hat{x}_{k-i-l+2}),
\end{equation}
where $\hat{B}_{0,0} = 1$ and $\hat{B}_{k,0}=0$ for $k \geq 1$.

\textbf{Proof.}
Using the relation between the exponential and ordinary Bell polynomials \eqref{pobp2}, the formula \eqref{pebp4} changes to
\begin{equation*}
\frac{k!}{l!} \hat{B} _{k,l} \left( \frac{x_1}{1!}, \ldots, \frac{x_{k-l+1}}{(k-l+1)!} \right)  = \sum\limits_{i=1}^{k-l+1} \binom{k-1}{i-1} x_i \frac{(k-i)!}{(l-1)!} \hat{B}_{k-i,l-1} \left( \frac{x_1}{1!}, \ldots, \frac{x_{k-i-l+2}}{(k-i-l+2)!} \right) .
\end{equation*}
After rearranging the expression and multiplying $i$th term in the sum by ${\displaystyle \frac{i!}{i!}}$, we obtain
\begin{align*}
\hat{B} _{k,l} \biggl( \frac{x_1}{1!}, \ldots &, \frac{x_{k-l+1}}{(k-l+1)!} \biggr) =\\
 &= \sum\limits_{i=1}^{k-l+1} \frac{(k-1)!}{(i-1)! (k-i)!} \frac{l!}{k!}\frac{(k-i)!}{(l-1)!} i! \frac{x_i}{i!} \hat{B}_{k-i,l-1} \left( \frac{x_1}{1!}, \ldots, \frac{x_{k-i-l+2}}{(k-i-l+2)!} \right).
\end{align*}
Now we cancel all possible factors and factorials and we get
\begin{equation}
\hat{B}_{k,l} \left (\frac{x_1}{1!}, \ldots, \frac{x_{k-l+1}}{(k-l+1)!} \right) = \sum\limits_{i=1}^{k-l+1} \frac{i \cdot l}{k} \frac{x_i}{i!} \hat{B}_{k-i,l-1} \left (\frac{x_1}{1!}, \ldots, \frac{x_{k-i-l+2}}{(k-i-l+2)!} \right).
\end{equation}
Denoting $\displaystyle{\hat{x}_i = \frac{x_i}{i!}}$ for all $i = 1, \ldots, k - l + 1$, gives the formula \eqref{pobp4}.

The main result of the paper is formulated in the following theorem:

\textbf{Theorem 2}
\label{mainth}
Let $g$ and $f$ be real functions analytic near $t_0$ and $g(t_0)$ respectively, and let $h$ be the composition $ h(t) = (f \circ g)(t)=f(g(t))$. Denote $\mathcal{D} \{ g(t) \}[t_0] = \{ G(k) \}_{k=0}^{\infty}$, $\mathcal{D} \{ f(t) \}[g(t_0)] = \{ F(k) \}_{k=0}^{\infty}$ and $\mathcal{D} \{ (f \circ g)(t) \}[t_0] = \{ H(k) \}_{k=0}^{\infty}$ the differential transformations of functions $g$, $f$ and $h$ at $t_0$, $g(t_0)$ and $t_0$ respectively. Then the numbers $H(k)$ in the sequence $\{ H(k) \}_{k=0}^{\infty}$ satisfy the relations $H(0) = F(0)$ and
\begin{equation}\label{fdb2}
H(k) = \sum_{l=1}^{k} F(l) \cdot \hat{B}_{k,l} \bigl( G(1), \ldots, G(k-l+1) \bigr) \ \text{ for } k \geq 1.
\end{equation}

\textbf{Proof.}
The assumption of analyticity of functions $g$, $f$ and $h$  guarantees that Theorem 1 is valid. Applying Definition 1 on formula \eqref{fdb}, we obtain
\begin{equation}
k! H(k) = \sum_{l=1}^{k} l! F(l) \cdot B_{k,l} \bigl( 1! G(1), \ldots, (k-l+1)! G(k-l+1) \bigr).
\end{equation}
If we divide both sides by $k!$ and apply formula \eqref{pobp2} in Lemma 4, we have
\begin{equation}
H(k) = \sum_{l=1}^{k} \frac{l!}{k!} F(l) \cdot \frac{k!}{l!} \hat{B}_{k,l} \left( \frac{1! G(1)}{1!}, \ldots, \frac{(k-l+1)! G(k-l+1)}{(k-l+1)!} \right)
\end{equation}
for $k \geq 1$. Cancelling all possible factorials gives the result.

\section{Applications}
\label{examples}
To show efficiency of the derived algorithm, we apply the differential transformation to two differential equations nonlinear with respect to the dependent variable $u(t)$.

\textbf{Example 1}
Let us consider equation
\begin{equation}\label{ex1.1}
u'(t) = u(t) - t + \ln (u(t))
\end{equation}
with initial condition
\begin{equation}\label{ex1.2}
u(0) = 1.
\end{equation}
Here we denote $h(t) = f(g(t))$, where $g(t) = u(t)$ and $f(x) = \ln (x)$. We are looking for an analytic solution in a neighbourhood of $t_0=0$. Then the righthand side of \eqref{ex1.1} is analytic too. The differential transformation turns the equation \eqref{ex1.1} into
\begin{equation}\label{ex1.3}
(k+1) U(k+1)[0] = U(k)[0] - \delta (k-1) + H(k)[0],
\end{equation}
with transformed initial condition $U(0)[0]=1$.

To find the coefficients $H(k)[0]$, we use Theorem 2. First of all, we recall Lemma 1, formula $v)$, to see that the differential transformation of $f(x) = \ln (x)$ at $x_0 = u(0) = 1$ is $F(0)[1]=0$ and $F(k)[1] = {\displaystyle \frac{(-1)^{k+1}}{k}}$ for $k \geq 1$. Theorem 2 gives $H(0)[0] = F(0)[1] = 0$ and\\
${\displaystyle H(k)[0] = \sum\limits_{l=1}^k F(l)[1] \cdot \hat{B}_{k,l} \bigl( U(1)[0], \ldots, U(k-l+1)[0] \bigr) \ \text{ for } k \geq 1.}$
We calculate
\begin{eqnarray*}
H(1)[0] &=& F(1)[1] \cdot \hat{B}_{1,1} \bigl( U(1)[0] \bigr) = U(1)[0],  \\
H(2)[0] &=&  \sum\limits_{l=1}^2 F(l)[1] \cdot \hat{B}_{2,l} \bigl( U(1)[0], U(2)[0] \bigr) = F(1)[1] \hat{B}_{2,1} \bigl( U(1)[0], U(2)[0] \bigr) +\\
&+& F(2)[1] \hat{B}_{2,2} \bigl( U(1)[0] \bigr)= U(2)[0] + \frac{(-1)}{2} \bigl( U(1)[0] \bigr)^2 = U(2)[0] - \frac{1}{2} \bigl( U(1)[0] \bigr)^2,\\
H(3)[0] &=& \sum\limits_{l=1}^3 F(l)[1] \cdot \hat{B}_{3,l} \bigl( U(1)[0], U(2)[0], U(3)[0] \bigr) \\
 &=& F(1)[1] \hat{B}_{3,1} \bigl( U(1)[0], U(2)[0], U(3)[0] \bigr)+ F(2)[1] \hat{B}_{3,2} \bigl( U(1)[0], U(2)[0] \bigr) +\\
 &+& F(3)[1] \hat{B}_{3,3} \bigl( U(1)[0] \bigr)= U(3)[0] - \frac{1}{2} 2 U(1)[0] U(2)[0] + \frac{1}{3} \bigl( U(1)[0] \bigr)^3.
\end{eqnarray*}
Now we are prepared to calculate the coefficients $U(k)[0]$, $k=1, 2, \ldots$.
\begin{eqnarray*}
U(1)[0] &=& U(0)[0] - \delta (-1) + H(0)[0] = 1,  \\
U(2)[0] &=&  \frac{1}{2} = \bigl( U(1)[0] - \delta (0) + H(1)[0] \bigr)= \frac{1}{2} ( 1 - 1 + 1) = \frac{1}{2},\\
U(3)[0] &=& \frac{1}{3} \bigl( U(2)[0] - \delta (1) + H(2)[0] \bigr) \!= \! \frac{1}{3} \! \left( \! \frac{1}{2} +U(2)[0] - \frac{1}{2} \bigl( U(1)[0] \bigr)^2 \! \right) \! \! = \! \frac{1}{3 \cdot 2} \! = \! \frac{1}{3!},\\
U(4)[0] &=& \frac{1}{4} \bigl( U(3)[0] - \delta (2) + H(3)[0] \bigr)= \frac{1}{4} \left( \frac{1}{3!} +U(3)[0] - \frac{1}{2} 2 U(1)[0] U(2)[0] + \right. \\
 &+& \left. \frac{1}{3} \bigl( U(1)[0] \bigr)^3 \right)= \ldots = \frac{1}{4!}, \\
 &\vdots&
\end{eqnarray*}
It is not difficult to verify that the $k$th coefficient $U(k)[0]$ has the value $\displaystyle{\frac{1}{k!}}$. Using the inverse differential transformation (Definition 2) we obtain
\begin{equation}
u(t) = \mathcal D^{-1} \Bigl\{ \ \{ U(k) [0] \}_{k=0}^{\infty} \ \Bigr\} [0] = \sum\limits_{k=0}^{\infty} \frac{1}{k!} t^k = \e^t,
\end{equation}
which is the exact solution of the initial value problem \eqref{ex1.1}, \eqref{ex1.2}.

\textbf{Example 2}
Let us solve the following equation
\begin{equation}\label{ex2.1}
u'(t) = 2 \sqrt{1- u^2 (t)}
\end{equation}
with initial condition
\begin{equation}\label{ex2.2}
u(0) = 0.
\end{equation}
We denote $h(t) = f(g(t))$, where $g(t) = u(t)$ and $f(x) = \sqrt{1-x^2}$. We want to find an analytic solution in a neighbourhood of $t_0=0$. The righthand side of \eqref{ex2.1} is analytic near $t_0$. The differential transformation of the equation \eqref{ex2.1} is
\begin{equation}\label{ex2.3}
(k+1) U(k+1)[0] = 2 H(k)[0],
\end{equation}
with transformed initial condition $U(0)[0]=0$.

We use Theorem 2 to find the coefficients $H(k)[0]$. First of all, we use Lemma 1, formula $vi)$ to derive the differential transformation of $f(x)$ at $x_0 = u(0) = 0$. Since ${\displaystyle f(x) = \sqrt{1-x^2} = (1+ (-x^2))^{\frac{1}{2}}}$, Taylor series of $f(x)$ at $0$ is ${\displaystyle f(x) = \sum\limits_{k=0}^{\infty} \binom{\frac{1}{2}}{k} (-x^2)^k}$. Consequently, $F(2k)[0]= {\displaystyle \binom{\frac{1}{2}}{k} (-1)^k}$ and $F(2k+1)[0] = 0$ for $k \in \mathbb{N}_0$. Theorem 2 gives $H(0)[0] = F(0)[0] = 1$ and
\begin{equation}\label{ex2.4}
H(k)[0] = \sum\limits_{l=1}^k F(l)[0] \cdot \hat{B}_{k,l} \bigl( U(1)[0], \ldots, U(k-l+1)[0] \bigr) \ \text{ for } k \geq 1.
\end{equation}
Now we take \eqref{ex2.3} and \eqref{ex2.4} in turns to obtain the coefficients.
\begin{eqnarray*}
U(1)[0] &=& \frac{2}{1} H(0)[0] = 2, \\
H(1)[0] &=& F(1)[0] \cdot \hat{B}_{1,1} \bigl( U(1)[0] \bigr) = 0 \cdot U(1)[0] = 0,  \\
U(2)[0] &=& \frac{2}{2} H(1)[0] = 0, \\
H(2)[0] &=&  \sum\limits_{l=1}^2 F(l)[0] \cdot \hat{B}_{2,l} \bigl( U(1)[0], U(2)[0] \bigr) = F(1)[0] \hat{B}_{2,1} \bigl( U(1)[0], U(2)[0] \bigr) +\\
&+& F(2)[0] \hat{B}_{2,2} \bigl( U(1)[0] \bigr)= 0 \cdot U(2)[0] + \frac{1}{2} (-1) \bigl( U(1)[0] \bigr)^2= - \frac{1}{2} 2^2,\\
U(3)[0] &=& \frac{2}{3} H(2)[0] = \frac{2}{3} \left( - \frac{1}{2} \right) 2^2 = - \frac{2^3}{3!}, \\
H(3)[0] &=& \sum\limits_{l=1}^3 F(l)[0] \cdot \hat{B}_{3,l} \bigl( U(1)[0], U(2)[0], U(3)[0] \bigr) \\
 &=& F(1)[0] \hat{B}_{3,1} \bigl( U(1)[0], U(2)[0], U(3)[0] \bigr)+ F(2)[0] \hat{B}_{3,2} \bigl( U(1)[0], U(2)[0] \bigr) +\\
 &+& F(3)[0] \hat{B}_{3,3} \bigl( U(1)[0] \bigr)= 0 - \frac{1}{2} 2 U(1)[0] U(2)[0] + 0 = 0, \\
U(4)[0] &=& \frac{2}{4} H(3)[0] = 0,\\
H(4)[0] &=& \sum\limits_{l=1}^4 F(l)[0] \cdot \hat{B}_{4,l} \bigl( U(1)[0], U(2)[0], U(3)[0], U(4)[0] \bigr) \\
 &=& F(1)[0] \hat{B}_{4,1} \bigl( U(1)[0], U(2)[0], U(3)[0], U(4)[0] \bigr)+ \\
 &+& F(2)[0] \hat{B}_{4,2} \bigl( U(1)[0], U(2)[0], U(3)[0] \bigr) + F(3)[0] \hat{B}_{4,3} \bigl( U(1)[0], U(2)[0] \bigr) +\\
 &+& F(4)[0] \hat{B}_{4,4} \bigl( U(1)[0] \bigr) \! = 0 - \frac{1}{2} \bigl( 2 U(1)[0] U(3)[0] + (U(2)[0])^2 \bigr) + 0 + \binom{\frac{1}{2}}{2} \bigl( U(1)[0] \bigr)^4 \\
 &=& - \frac{1}{2}  \bigl( 2 \cdot 2 \cdot \frac{-2^3}{3!} + 0 \bigr) + \frac{1}{2} \cdot \left( - \frac{1}{2} \right) \cdot \frac{1}{2!} \cdot 2^4 = \frac{2^4}{6} - \frac{2^4}{8} = \frac{2^4}{24} = \frac{2^4}{4!}, \\
U(5)[0] &=& \frac{2}{5} H(4)[0] = \frac{2}{5} \cdot \frac{2^4}{4!} = \frac{2^5}{5!}, \\
 &\vdots&
\end{eqnarray*}
Also in this case, it is possible to see the pattern for the $k$th coefficient $U(k)[0]$, which is $U(k)[0] =0$ if $k$ is even and  ${\displaystyle U(k)[0] = (-1)^{\frac{k-1}{2}}\frac{2^k}{k!}}$ if $k$ is odd. The inverse differential transform gives
\begin{equation*}
u(t) \! = \! \mathcal D^{-1} \Bigl\{ \! \{ U(k) [0] \}_{k=0}^{\infty} \Bigr\} \! [0] \! = \! \! \sum\limits_{k=0}^{\infty} \! \frac{(-1)^k \cdot 2^{2k+1}}{(2k+1)!} t^{2k+1} \! = \! \! \sum\limits_{k=0}^{\infty} \! \frac{(-1)^k}{(2k+1)!} (2t)^{2k+1} \! = \! \sin 2t.
\end{equation*}
Indeed, $u(t) = \sin 2t$ is the exact unique solution of the given initial value problem \eqref{ex2.1}, \eqref{ex2.2}.

\textbf{Remark}\label{remark6}
 In recent works \cite{petro1}, \cite{petro2}, a more general discretization technique involving Hilbert spaces is introduced. Although the approach is different, the transformation rules are identical to the differential transformation formulas.

We can summarize the main advantages of the presented approach as follows:
\begin{itemize}
\item
Using the presented algorithm, we are able to obtain approximate solution of the initial value problem. However, in some cases, there is the possibility to identify the unique solution in closed form.
\item
We do not need initial guess approximation and symbolic computation of multiple integrals or derivatives, hence less calculations are demanded compared to other popular semi-analytical methods (the variational iteration method, the homotopy perturbation method, the homotopy analysis method, the Adomian decomposition method).
\item
There is no need for numerical integration or differentiation either. Only arithmetical operations are used.
\item
In comparison to any purely numerical method, a specific advantage of this technique is that the approximate solution is always a function analytic near $t_0$.
\item
The algorithm is recurrent, so we use values computed in previous steps. Suitable arrangements can be made to reduce the necessary computational work. This fact was demonstrated in Example 2.
\end{itemize}

\section{Conclusion}
Results of this paper are based on combinatorial properties of Bell polynomials applied in the differential transformation theory. Modified version of the Fa\`{a} di Bruno's formula with partial ordinary Bell polynomials was proved. Applicability of the algorithm was demonstrated on two particular examples of the initial value problem for differential equations with nonlinearity containing the unknown function. The algorithm can be generalized to other types of problems, e.g. boundary value problems.

\section{Acknowledgements}
\label{sec6}
The author was supported by the Grant FEKT-S-17-4225 of Faculty of Electrical Engineering and Communication, Brno University of Technology. This support is gratefully acknowledged.




\section*{Current address}

\vspace{1.2ex}
\usebox{\authors}

\end{aplart}
\end{document}